\documentclass{amsart}
\usepackage{amsmath,amscd,amsthm,amsfonts,amssymb}

\theoremstyle{plain}
\newtheorem{theorem}                 {Theorem}      [section]
\newtheorem{proposition}  [theorem]  {Proposition}

\theoremstyle{definition}
\newtheorem{example}      [theorem]  {Example}

\newtheorem{definition}   [theorem]  {Definition}

\numberwithin{equation}{section}

\def \theo-intro#1#2 {\vskip .25cm\noindent{\bf Theorem #1\ }{\it #2}}

\def \rn{\mathbb R}
\def \cn{\mathbb C}

\def \F{\mathcal F}
\def \H{\mathcal H}

\def \V{\mathcal V}

\def \ip #1#2{\langle #1,#2 \rangle}

\def \lb#1#2{[#1,#2]}

\def \g{\mathfrak{g}}

\def \nab#1#2{\hbox{$\nabla$\kern -.3em\lower 1.0 ex
    \hbox{$#1$}\kern -.1 em {$#2$}}}

\begin{document}
\baselineskip 22pt \larger

\allowdisplaybreaks

\title
{Holomorphic harmonic morphisms from four-dimensional non-Einstein manifolds}

\author{Sigmundur Gudmundsson}



\dedicatory
{version 1.008 - 21 May 2014}

\address
{Department of Mathematics, Faculty of Science, Lund University,
Box 118, S-221 00 Lund, Sweden}

\email
{Sigmundur.Gudmundsson@math.lu.se}

\begin{abstract}
We construct 4-dimensional Riemannian Lie groups carrying left-invariant
conformal foliations with minimal leaves of codimension 2. We show that
these foliations are holomorphic with respect to an (integrable) Hermitian
structure which is {\it not} K\" ahler.  We then prove that the Riemannian Lie
groups constructed are {\it not} Einstein manifolds.  This answers an important
open question in the theory of complex-valued harmonic morphisms from
Riemannian 4-manifolds.
\end{abstract}

\maketitle

\section{Introduction}

More than twenty years ago, J. C. Wood proved, in \cite{Woo}, that any submersive
harmonic morphism from an orientable 4-dimensional Einstein manifold $M$
to a Riemann surface, or a conformal foliation of $M$ by minimal surfaces,
determines an (integrable) Hermitian structure with respect to which it
is holomorphic. Ever since it has been an open question whether
holomorphicity in the above situation forces the 4-manifold $M$ to be Einstein.

In this paper we construct two 3-dimensional families of 4-dimensional
Riemannian Lie groups carrying left-invariant conformal foliations with
minimal leaves of codimension 2. We show that these foliations are holomorphic
with respect to an (integrable) Hermitian structure which is {\it not} K\" ahler.
We then prove that the Riemannian Lie groups constructed are {\it not}
Einstein manifolds.  This gives a definite answer to the above mentioned
open question.

For the general theory of harmonic morphisms between Riemannian
manifolds we refer to the excellent book \cite{Bai-Woo-book}
and the regularly updated on-line bibliography \cite{Gud-bib}.

\section{Harmonic morphisms and minimal conformal foliations}

Let $M$ and $N$ be two manifolds of dimensions $m$ and $n$,
respectively. A Riemannian metric $g$ on $M$ gives rise to the
notion of a {\it Laplacian} on $(M,g)$ and real-valued {\it
harmonic functions} $f:(M,g)\to\rn$. This can be generalized to
the concept of {\it harmonic maps} $\phi:(M,g)\to (N,h)$ between
Riemannian manifolds, which are solutions to a semi-linear system
of partial differential equations, see \cite{Bai-Woo-book}.

\begin{definition}
  A map $\phi:(M,g)\to (N,h)$ between Riemannian manifolds is
  called a {\it harmonic morphism} if, for any harmonic function
  $f:U\to\rn$ defined on an open subset $U$ of $N$ with $\phi^{-1}(U)$
non-empty,
  $f\circ\phi:\phi^{-1}(U)\to\rn$ is a harmonic function.
\end{definition}

The following characterization of harmonic morphisms between
Riemannian manifolds is due to Fuglede and T. Ishihara.  For the
definition of horizontal (weak) conformality we refer to
\cite{Bai-Woo-book}.

\begin{theorem}\cite{Fug-1,T-Ish}
  A map $\phi:(M,g)\to (N,h)$ between Riemannian manifolds is a
  harmonic morphism if and only if it is a horizontally (weakly)
  conformal harmonic map.
\end{theorem}

Let $(M,g)$ be a Riemannian manifold, $\V$ be an involutive
distribution on $M$ and denote by $\H$ its orthogonal
complement distribution on $M$.
As customary, we also use $\V$ and $\H$ to denote the
orthogonal projections onto the corresponding subbundles of $TM$
and denote by $\F$ the foliation tangent to
$\V$. The second fundamental form for $\V$ is given by
$$B^\V(U,V)=\frac 12\H(\nabla_UV+\nabla_VU)\qquad(U,V\in\V),$$
while the second fundamental form for $\H$ satisfies
$$B^\H(X,Y)=\frac{1}{2}\V(\nabla_XY+\nabla_YX)\qquad(X,Y\in\H).$$
The foliation $\F$ tangent to $\V$ is said to be {\it conformal} if there is a
vector field $V\in \V$ such that $$B^\H=g\otimes V,$$ and
$\F$ is said to be {\it Riemannian} if $V=0$.
Furthermore, $\F$ is said to be {\it minimal} if $\text{trace}\ B^\V=0$ and
{\it totally geodesic} if $B^\V=0$. This is equivalent to the
leaves of $\F$ being minimal and totally geodesic submanifolds
of $M$, respectively.

It is easy to see that the fibres of a horizontally conformal
map (resp.\ Riemannian submersion) give rise to a conformal foliation
(resp.\ Riemannian foliation). Conversely, the leaves of any
conformal foliation (resp.\ Riemannian foliation) are
locally the fibres of a horizontally conformal map
(resp.\ Riemannian submersion), see \cite{Bai-Woo-book}.

The next result of Baird and Eells gives the theory of
harmonic morphisms, with values in a surface,
a strong geometric flavour.

\begin{theorem}\cite{Bai-Eel}\label{theo:B-E}
Let $\phi:(M^m,g)\to (N^2,h)$ be a horizontally conformal
submersion from a Riemannian manifold to a surface. Then $\phi$ is
harmonic if and only if $\phi$ has minimal fibres.
\end{theorem}

\section{4-dimensional Lie groups}

Let $G$ be a 4-dimensional Lie group equipped with a left-invariant
Riemannian metric.  Let $\g$ be the Lie algebra of $G$ and
$\{X,Y,Z,W\}$ be an orthonormal basis for $\g$.  Let $Z,W\in\g$
generate a 2-dimensional left-invariant and integrable distribution $\V$
on $G$ which is conformal and with minimal leaves.  We denote
by $\H$ the horizontal distribution, orthogonal to $\V$, generated by
$X,Y\in\g$.  Then it is easily seen that the basis $\{ X,Y,Z,W\}$ can
be chosen so that the Lie bracket relations for $\g$ are of the form
\begin{eqnarray*}
\lb WZ&=&\lambda W,\\
\lb ZX&=&\alpha X +\beta Y+z_1 Z+w_1 W,\\
\lb ZY&=&-\beta X+\alpha Y+z_2 Z+w_2 W,\\
\lb WX&=&     a X     +b Y+z_3 Z-z_1W,\\
\lb WY&=&    -b X     +a Y+z_4 Z-z_2W,\\
\lb YX&=&     r X         +\theta_1 Z+\theta_2 W
\end{eqnarray*}
with real structure constants.  For later reference we state the following
easy result describing the geometry of the situation.

\begin{proposition}\label{prop-geometry}
Let $G$ be a 4-dimensional Lie group and $\{X,Y,Z,W\}$ be an
orthonormal basis for its Lie algebra as  above.  Then
\begin{enumerate}
\item[(i)] $\F$ is {\it totally geodesic} if and only if
  $z_1=z_2=z_3+w_1=z_4+w_2=0$,
\item[(ii)] $\F$ is {\it Riemannian} if and only if $\alpha=a=0$, and
\item[(iii)] $\H$ is {\it integrable} if and only if $\theta_1=\theta_2=0$.
\end{enumerate}
\end{proposition}

On the Riemannian Lie group $G$ there exist, up to sign,
exactly two invariant almost Hermitian structure $J_1$ and $J_2$
which are adapted to the orthogonal decomposition $\g=\V\oplus\H$
of the Lie algebra $\g$.  They are determined by
$$J_1X=Y,\ J_1Y=-X,\ J_1Z=W,\ J_1W=-Z,$$
$$J_2X=Y,\ J_2Y=-X,\ J_2W=Z,\ J_2Z=-W.$$
An elementary calculation involving the Nijenhuis tensor shows that
$J_1$ is integrable if and only if
\begin{equation}
2z_1-z_4-w_2=2z_2+z_3+w_1=0.
\end{equation}
In this case the Lie bracket relations for $\g$ take the form
\begin{eqnarray*}
\lb WZ&=&\lambda W,\\
\lb ZX&=&\alpha X +\beta Y+z_1 Z-(2z_2+z_3)W,\\
\lb ZY&=&-\beta X+\alpha Y+z_2 Z+(2z_1-z_4)W,\\
\lb WX&=&     a X     +b Y+z_3 Z-z_1W,\\
\lb WY&=&    -b X     +a Y+z_4 Z-z_2W,\\
\lb YX&=&     r X         +\theta_1 Z+\theta_2 W.
\end{eqnarray*}

\begin{example}
With the non-vanishing coefficients $z_2,\theta_1,\theta_2$ we
obtain the following 3-dimensional family of Lie algebras
\begin{eqnarray*}
\lb ZX&=&-2z_2 W,\\
\lb ZY&=&  z_2 Z,\\
\lb WY&=& -z_2 W,\\
\lb YX&=& 2z_2 X+\theta_1 Z+\theta_2 W.
\end{eqnarray*}
These are the special cases $\g_{7}(z_2,-2z_2,0,\theta_1,\theta_2)$ of
Example 5.3 in \cite{Gud-Sve-6}.  It should be noted that the horizontal 
distribution $\H$ is not integrable and the leaves of the vertical
foliation $\V$ are not totally geodesic.

We will show that none of the corresponding Riemannian Lie groups
are Einstein manifolds. A standard calculation involving the Koszul fomula
$$2\ip{\nab XY}Z=\ip{\lb ZX}Y+\ip{\lb ZY}X+\ip Z{\lb XY}$$
shows that the Levi-Civita connection satisfies the following relations
$$\nab XX=2z_2Y,\ \ \nab XY=-2z_2X-\frac 12\theta_1Z-\frac 12\theta_2W,$$
$$\nab XZ=\frac 12\theta_1Y+z_2W, \ \ \nab XW=\frac 12\theta_2Y-z_2Z,$$
$$\nab YX=\frac 12\theta_1Z+\frac 12\theta_2W,\ \ \nab YY=0,$$
$$\nab YZ=-\frac 12\theta_1X,\ \ \nab YW=-\frac 12\theta_2X,$$
$$\nab ZX=\frac 12\theta_1Y-z_2W,\ \ \nab ZY=-\frac 12\theta_1X+z_2Z,$$
$$\nab ZZ=-z_2Y,\ \ \nab ZW=z_2X,$$
$$\nab WX=\frac 12\theta_2Y-z_2Z,\ \ \nab WY=-\frac 12\theta_2X-z_2W,$$
$$\nab WZ=z_2X,\ \ \nab WW=z_2Y.$$
This means that the Hermitian structure $J_1$ is not K\" ahler,
since $$(\nab X{J_1})(X)=-\frac 12(\theta_1Z+\theta_2W).$$
Employing the definition for the sectional curvature
$$\ip{R(X,Y)Y}X=\ip{\nab X{\nab YY}-\nab Y{\nab XY}-\nab{[X,Y]}Y}X$$
we then obtain the following useful equalities
$$\ip{R(X,Y)Y}X=-\frac 34(\theta_1^2+\theta_2^2)-4z_2^2,\ \ 
\ip{R(X,Z)Z}X=\frac 14\theta_1^2-z_2^2,$$
$$\ip{R(X,W)W}X=\frac 14\theta_2^2-z_2^2,\ \ 
\ip{R(Y,Z)Z}Y=\frac 14\theta_1^2-z_2^2,$$
$$\ip{R(Y,W)W}Y=\frac 14\theta_2^2-z_2^2,\ \ 
\ip{R(Z,W)W}Z=2z_2^2.$$
One immediate consequence is that
$$Ric (X,X)=-\frac 12(\theta_1^2+\theta_2^2)-6z_2^2,\ \ 
Ric(Z,Z)=\frac 12\theta_1^2.$$
showing that none of these Riemannian Lie groups is an Einstein manifold.
\end{example}

\section{The structures $J_1$ and $J_2$ are both integrable}

It is easily seen that $J_1$ and $J_2$ are {\it both} integrable if and only if
$$2z_1-z_4-w_2=2z_2+z_3+w_1=0,$$
$$2z_1+z_4+w_2=2z_2-z_3-w_1=0.$$
As a direct consequence, we see that in this case the foliation $\F$
is totally geodesic i.e.
$$z_1=z_2=z_3+w_1=z_4+w_2=0.$$
In this situation the Lie bracket relations take the following form
\begin{eqnarray*}
\lb WZ&=&\lambda W,\\
\lb ZX&=&\alpha X +\beta Y          -z_3 W,\\
\lb ZY&=&-\beta X+\alpha Y          -z_4 W,\\
\lb WX&=&     a X     +b Y+z_3 Z,\\
\lb WY&=&    -b X     +a Y+z_4 Z,\\
\lb YX&=&     r X    +\theta_1 Z+\theta_2 W.
\end{eqnarray*}

\begin{example}
For the non-vanishing coefficients $\alpha,\beta,\theta_2$ we have the
following 3-dimensional family of Lie algebras
\begin{eqnarray*}
\lb WZ&=&-2\alpha W,\\
\lb ZX&=&\alpha X +\beta Y,\\
\lb ZY&=&-\beta X+\alpha Y,\\
\lb YX&=&\theta_2 W.
\end{eqnarray*}
They are the special cases $\g_{3}(\alpha,\beta,0,0,\theta_2)$ of
Example 4.3 in \cite{Gud-Sve-6}.  For the corresponding Riemannian
Lie groups the Levi-Civita connection is given by
$$\nab XX=\alpha Z,\ \ \nab XY=-\frac 12\theta_2W,\ \
\nab XZ=-\alpha X,\ \ \nab XW=\frac 12\theta_2Y,$$
$$\nab YX=\frac 12\theta_2W,\ \ \nab YY=\alpha Z,\ \
\nab YZ=-\alpha Y,\ \ \nab YW=-\frac 12\theta_2 X,$$
$$\nab ZX=\beta Y,\ \ \nab ZY=-\beta X,\ \
\nab ZZ=0,\ \ \nab ZW=0,$$
$$\nab WX=\frac 12\theta_2Y,\ \ \nab WY=-\frac 12\theta_2X,\ \
\nab WZ=-2\alpha W,\ \ \nab WW=2\alpha Z.$$
Using these identities it is easily seen that
\begin{itemize}
\item[i)] the Hermitian structure $J_1$ is K\"ahler if and only if $\theta_2=-2\alpha$,
\item[ii)] the Hermitian structure $J_2$ is K\"ahler if and only if $\theta_2=2\alpha$.
\end{itemize}
This means that in all the cases that we are considering at least 
one of the Hermitian structures $J_1$ or $J_2$ is not K\" ahler. 
For the sectional curvatures have
$$\ip{R(X,Y)Y}X=-\alpha^2-\frac 34\theta_2^2,\ \ \ip{R(X,Z)Z}X=-\alpha^2,$$
$$\ip{R(X,W)W}X=\frac 14\theta_2^2-2\alpha^2,\ \ \ip{R(Y,Z)Z}Y=-\alpha^2,$$
$$\ip{R(Y,W)W}Y=\frac 14\theta_2^2-2\alpha^2,\ \ \ip{R(Z,W)W}Z=-4\alpha^2.$$
Finally, we yield 
$$Ric (X,X)=-\frac 12\theta_2^2-4\alpha^2,\ \ 
Ric(Z,Z)=-6\alpha^2,\ \ Ric(W,W)=\frac 12\theta_2^2-8\alpha^2$$
telling us that if $4\alpha^2\neq \theta_2^2$ then our Riemannian 
Lie group is not an Einstein manifold.
\end{example}

\section{Acknowledgements}
The author is grateful to John $\cn.$ Wood for useful discussions
on this work at the {\it Differential Geometry Workshop} held at
Lund in May 2014.

\end{document}